\documentclass[11pt]{article}
\usepackage{graphicx,epsfig,psfrag}
\usepackage{amssymb,amsmath,amscd,amsthm}
\usepackage{fullpage}
 \newtheorem{thm}{Theorem}  
  \newtheorem{cor}{Corollary}

 \def\illustration #1 by #2 (#3)
(#4){\leavevmode
  \vbox to #2{
   \hrule width #1 height 0pt depth 0pt
   \vfill
   \special{illustration #3 scaled #4}}} 
\begin{document} 
\title{Minimal capacity points and the lowest eigenfunctions.} 
\author{ Mark Levi and Jia Pan\thanks{Research of both authors was partially supported  by 
the NSF grant DMS-9704554 } 
  } 
\bibliographystyle{plain} 
\maketitle

We introduce the concept of the  a point of minimal capacity of the domain, and observe   a connection between  this point and  the lowest eigenfunction of a Laplacian on a domain, in one special case.  

\section{A physical motivation and precise definitions}
\subsection{The ``warmest" point of a domain.}
For a given a domain in $D\subset {\mathbb R}  ^n$, let  $ m(D)\in D$ be {\it  a point of maximum of   the lowest eigenfunction of the Laplacian with zero Dirichlet boundary conditions\footnote{Such a point need not be unique for non-convex domains.}.}  A point $ m(D) $ is in a certain sense the warmest point in $D$: if the domain, viewed as a heat-conducting medium, starts with a positive temperature distribution and is cooled by   the maintenance of zero temperature on the boundary, 
  then after a long time a (local) maximum of the temperature will approach a maximum  of the lower eigenfunction.
Indeed, the solution $ u $  of the heat equation 
$$
u_t=\Delta u, \ u(0,x)=u_0(x), \ u(t,x)|_{x\in\partial D}=0 
$$
on a     domain $D$ in $ {\mathbb R} ^n$ is   given in terms of eigenfunctions $ v_n $ of the Laplacian with the Dirichlet boundary conditions: 
$$
u (x,t) = \sum_{n=0}^{\infty}a_ne^{ \lambda _nt} v_n(x),
$$
with $ 0> \lambda _0> \lambda _1>\cdots  $.     
The leading mode 
$ u(x,t)= e^{ \lambda_{0} t}v_0(x) $ becomes dominant for $t$ large (we assume $ a_0>0 $ and normalize to $ a_0=1$).  Hence   the warmest point   indeed approaches a maximum of $ u_0 $, as claimed. For the background on the Dirichlet problem we refer to the classical text \cite{courant-hilbert}. 
 
 \begin{figure}[thb]
         \psfrag{m}[c]{\small $m$}
 \psfrag{u1}[c]{\small $u=1$}
  \psfrag{u0}[c]{\small $u=0$}
   \psfrag{principal eigenfunction}[c]{\small principal eigenfunction}
    \psfrag{u}[c]{\small $u$}
    \psfrag{The open problem in 2D: c<m?}[c]{\small Open problem: is $c<m$?}
       \psfrag{D}[c]{\small $D$}
   \center{ \includegraphics[scale=0.5]{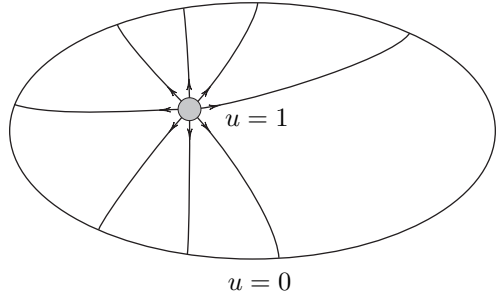}}
    \caption{The minimal capacity point $m$ minimizes the flux from the small disk around $m$ to $ \partial D$; a more precise definition is in the text.}
    \label{fig:2D}
\end{figure}

\subsection{The minimal capacity point in 2D.} 
Along with  the maximizer $ m (D) $  of the lowest eigenvalue we consider another special point which we will call  {\it  the minimal capacity point}  of the domain $D$. In a sense made precise shortly, this is the   ``best insulated from the boundary" point; we define it as follows. 
Let $ S_{ \varepsilon }(x)\subset D $ be the sphere of radius $\varepsilon$ centered at $ x\in D $, Figure~\ref{fig:2D}. Let us prescribe the boundary conditions  
$$
	u _{S_{ \varepsilon }(x)}= 1, \ \ {\rm and} \ \     u_{\partial D}=0 , 
$$ 
and 
let $u$ be  the solution of the resulting Dirichlet problem on the cored domain $D \setminus {\rm int}( S_{ \varepsilon }(x)) $.   If we interpret $u$ as the temperature, heat will flow from the warm sphere $ S_{ \varepsilon }(x) $ to the cold
boundary $ \partial D$, with the heat flux given by 
\begin{equation} 
     F_{ \varepsilon  }(x)=-\int_{S_{ \varepsilon }(x)} \nabla u \cdot {\bf n}\; dS.
          \label{eq:flux} 
\end{equation} 
This flux depends on the choice of $x$ and on $\varepsilon$. Let $ x_\varepsilon  $ be a minimizer of $ F_{ \varepsilon  } $
(such a minimizer may not be unique, for instance for a dumbbell--shaped domain with a thin neck).    If 
$$
	c(D) \buildrel{{\rm def}}\over{=}  \lim_{\varepsilon \rightarrow 0 }x_ \varepsilon 
$$ 
exists,  we will call this limit the  {\it   minimal capacity point}\footnote{No claim is made as to uniqueness of $ c(D) $; in fact, one can construct examples with multiple such points. In the case we consider, however, $ c $ is unique.}. Intuitively, $ c(D) $ is the point best insulated from the boundary, since it minimizes the flux of heat. 
\vskip 0.3 in 
An alternative physical interpretation:   $ c(D) $ is the {\it  point of least capacity of the capacitor whose electrodes are $ \partial D$ and $S_{x}( \varepsilon ) $}, in the limit of $ \varepsilon \rightarrow 0 $. Indeed,  if we interpret $ \nabla u $ as the electrostatic field in the vacuum (which we can do since $ \Delta u = 0 $), then the flux  $F_{ \varepsilon  }(x)$ of this field  is, by Gauss's law, and up to a scaling factor, the amount $q$ of electrostatic charge on the spherical electrode  $S_{x}( \varepsilon ) $. But the potential difference $V$ between the electrodes $ S_ \varepsilon (x) $ and $ \partial D$ is $ V=1 $, and 
by the definition of the capacitance ($C=q/V = F_ \varepsilon (x)/1 = F_ \varepsilon (x)$) we conclude that   $ F_ \varepsilon (x)$ is precisely  the capacity of the capacitor in question. 
\vskip 0.3 in

There is yet one more interpretation of $ F_ \varepsilon (x) $ defined by Eq. (\ref{eq:flux}): it is   the potential energy (up to a factor of $\frac{1}{2}$) stored in an elastic membrane subject to the boundary conditions mentioned in the definition of the least capacity point (for small displacements for which the nonlinear equation of the minimal surface can be replaced by $ \nabla u =0 $).  One can imagine  placing a heavy slippery disk $ S_ \varepsilon $  on the horizontal membrane. The disk   will slide to the position of least potential energy   (\ref{eq:flux}), i.e. to the vicinity of the least capacity point.   
 \vskip 0.3 in  
To summarize this and the preceding sections, we mentioned    two physically reasonable definitions of the ``warmest point" in $D$. The vague intuitive connection between these two definitions suggests a more precise mathematical relationship. The goal of this note is to explore this relationship on a simple example (Section \ref{mr}), and to state an open 
problem (Section \ref{op}). 
\vskip 0.3 in 
We mention in this connection   another problem suggested by Walter Craig and studied in depth by Jochen Denzler \cite{denzler1} \cite{denzler}. This problem addresses the question  ``where to place a window of given area to minimize heat loss?" The problem   reduces to 
 minimizing the principal eigenvalue of a mixed Neumann-Dirichlet problem.  

 \subsection{A characterization of the least capacity point. }
 Before discussing  the main result we give yet one more  characterization of the  minimal capacity point $ c(D)  $  in terms of Green's function of the domain.  
  \begin{thm} Let $D\subset {\mathbb R}  ^n $ be a domain for which any Dirichlet problem has a solution. 
Given  any $ p \in D$ let  $ v_{p}(z) $ be the harmonic function of $z$ satisfying Dirichlet boundary conditions    $ v_p= -\ln  \frac{1}{ |z-p|} $ on $ {\partial D} $  if $ n=2 $ and $ v_p= 
 - \frac{1}{|z-p|^{n-1}}  $ on $ {\partial D} $ if $ n>2 $. The least capacity point $ c(D)$ maximizes the function 
 $ v_z(z) $. 
 \end{thm}  
 \noindent{\bf Remark. \,}  The definition of   $ v_{p}(z) $ is motivated by the desire to make Green's function
 $$
 	  v_{p}(z) + \ln  \frac{1}{ |z-p|}
 $$ 
(for $ n =2$)  vanish on $ \partial D$. 
\vskip 0.3 in 

 \noindent{\bf Proof of Theorem 1.} We wish to construct the harmonic function involved in the definition of $ c(D)$. As a candidate, we take  Green's function\footnote{we consider the case $ n=2$, leaving out the obvious changes required for $ n > 2 $.} 
\begin{equation} 
    u (z) = k^{-1} \biggl( \ln |z-p|^{-1} + v_{p}(z) \biggr), \label{eq:defv}
\end{equation} 
where   $ p\in D$ and where $k$ will be chosen so that the average 
\begin{equation} 
     \overline{u}_{|z-p|=\varepsilon }=1. \label{eq:normailzation}
\end{equation} 
By the definition of $ v_p $ we have  $u_{\partial D}=0$. 
 
Averaging Eq. (\ref{eq:defv}) over the sphere $ S_ \varepsilon (p) $  and using the fact that   $v_p(z)$ is harmonic, we 
obtain the condition on $k$ which guarantees   (\ref{eq:normailzation}): 
 $$1= k ^{-1} (\ln1/ \varepsilon + v_{p}(p)) ,$$ 
or
$$ k= (\ln1/ \varepsilon + v_{p}(p) ) ^{-1}.$$
 Thus the harmonic function satisfying $ u_{\partial D}=0 $ and $ u_{|z-p|= \varepsilon }=1$ 
 is given by 
$$
u(z)= \frac{ \ln |z-p|^{-1} + v_{p}(z)}{\ln 1/ \varepsilon + v_{p}(p)} , 
$$
to the leading order for small $ \varepsilon$. 
Computing the ``heat flux"  we get the contribution $ 2 \pi $ from the logarithmic term 
and zero from $ v_p(z) $ since the latter is a harmonic function in $D$; to the leading order we have
 $$F_ \varepsilon (p)=- \oint_{|z-p|= \varepsilon } \nabla u\cdot {\bf n}  ds = \frac{2\pi }{\ln 1/ \varepsilon + v_{p}(p)}.$$
This shows that for small $\varepsilon$, the minimizer $ x_ \varepsilon $ of  $F_ \varepsilon (p) $ is close to the maximizer of  $  v_{p}(p) $. 
In other words,    $c(D) \buildrel{\rm def}\over{=} \lim_{\varepsilon \rightarrow 0 }x_{ \varepsilon }  $ is the maximizer $ v_p(p)$. 
 \hfill $\diamondsuit$   
\subsection {The point of minimal capacity for an elliptic operator on an interval. } 
For the case of $ n = 1 $ the concepts described before become trivial: both points $ m(D) $ and $ c(D) $ are simply the midpoints 
of the interval $D$.   However, the question is still interesting for a more general elliptic operator 
\begin{equation} 
Lu= (a(x) u ^\prime ) ^\prime , \ a(x)>0 \label{eq:oned}
\end{equation} 
on $ C^2[0,1]$ with the homogeneous boundary conditions $ u(0)=u(1)=0 $. 

Since this   is not a special case of the above, we repeat the definitions of the points $m$ and $c$. Since the domain 
can be taken as the unit interval $ [0,1]$, the points will essentially depend only on the operator $L$, and we will write 
 $ c[L] $ and $ m[L]$ in slight break with the earlier notation. 

\begin{figure}[thb]
         \psfrag{c}[c]{\small $c$}
     \psfrag{C}[c]{\small $s$}
              \psfrag{x}[c]{\small $x$}
                       \psfrag{u}[c]{\small $u$}
                       \psfrag{f}[c]{\small $F(c)$}
                       \psfrag{F}[c]{\small $F(s)$}
                       \psfrag{u-}[c]{\small $u_-(x;s)$}
                       \psfrag{u+}[c]{\small $u_+(x;s)$}
                       \psfrag{1}[c]{\small $1$}
  \center{   \includegraphics[scale=0.5]{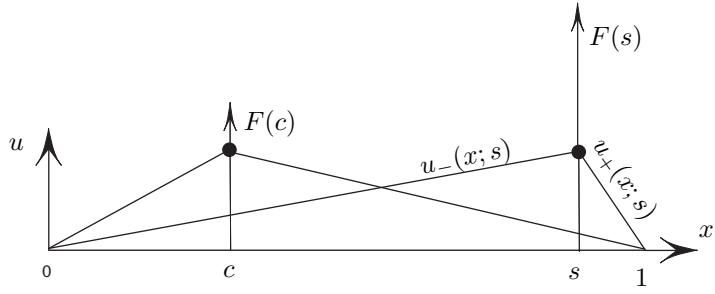} } 
    \caption{$F(s)>F(c) $: the flux from $c$ is smallest, i.e. $c$ is the best insulated point. As it turns out, $c$ also bisects electrical resistance of the segment if $ a(x) $ is interpreted as the local conductivity.   }
     \label{fig:string1}
\end{figure}

\paragraph {Definition of $ c[L]$.}  Let us introduce Green's function: for any $ 0<s<1$ (Figure~\ref{fig:string1}) we consider the solution 
$ u_-(x;s) $  of $ Lu=0  $ on  $  x\in[0, s]$ with $ u_-(0;s)=0 $, $u_-(s;s)=1$. Similarly we define $ u_+$  as the solution of $ Lu=0  $  on $ [s,1]$ with $ u_+(s;s)=1 $, $u_+(1;s)=0$.  We consider the (one-dimensional) flux out of $ x=s$, the analog of   (\ref{eq:flux}):  
\begin{equation} 	
	F(s)= a(s)(-u ^\prime  _+(x;s)+ u ^\prime _-(x;s))_{x=s}, 
	\label{eq:F}
\end{equation}   
where $ ^\prime $ denotes the $x$--derivative. Note that the signs in   (\ref{eq:F}) are chosen so that  the flux {\it  out} of $ x $ is counted with a positive sign, just like  in the higher dimensional case   (\ref{eq:flux}).    {\it  Finally, we define $c=c[L]$ as the minimizer of $F(s)$.}  It is clear that such a minimizer exists, since $F(s) \rightarrow \infty $ as $ s \rightarrow 0 $ or $ s \rightarrow 1 $. 

\vskip 0.3 in 
\noindent{\bf Remark 1. \,}  The least capacity point $ c=c[L] $  is alternatively  defined by 
\begin{equation} 
	\int_{0}^{c}  a ^{-1} (x) \;dx= \int_{c}^{1}    a ^{-1} (x)\;dx.
	\label{eq:conductivity}
\end{equation}   
	In particular,  the least capacity point is unique. 
	\vskip 0.3 in 
\noindent{\bf Proof.}    
From the definition of $ u_\pm $ as the solutions of $ Lu=0 $ with appropriate boundary conditions we obtain 
\begin{equation} 
	u_-(x;s)= \frac{\int_{0}^{x} a ^{-1}(t) dt}{\int_{0}^{s} a ^{-1}(t) dt},  \ \ \ 
	u_+(x;s)= \frac{\int_{x}^{1} a ^{-1}(t) dt}{\int_{s}^{1} a ^{-1}(t) dt}.
	\label{eq:uexplicit}
\end{equation}   
 Inserting these into   (\ref{eq:F}) we obtain 
 $$
 	F(s) =    \frac{1}{\int_{0}^{s} a ^{-1}(t) dt}+\frac{1}{\int_{s}^{1} a ^{-1}(t) dt} = \frac{1}{R(s)} + \frac{1}{R(1)-R(s)},  
 $$
 where $ R(s) = \int_{0}^{s} a ^{-1} (t)\;dt$. Now since $f(r) = \frac{1}{r} + \frac{1}{k-r} $ is maximized by 
 $ r = k/2 $, we conclude that $ F(s) $ is maximized by that value of $s$ which gives $ R(s) = \frac{1}{2} R(1) $, i.e. 
 by $ s=c $ satisfying   (\ref{eq:conductivity}).  \hfill $\diamondsuit$    
 
 \vskip 0.3 in 
 
 If $ a(x)$ is interpreted as the local electrical conductivity of a wire at $x$ (that is, the conductivity per unit length at $x$), then $ R(s) = \int_{0}^{s}a ^{-1} (t)\;dt$ is the resistance of the  piece of wire $ [0,s] $. 
 Note that    (\ref{eq:conductivity}) is intuitively plausible: it states that for $c$ to maximize the electrical resistance from itself to the   
two ends of the segment,  $c$ must {\it  bisect }  the resistance  of the entire segment. 
One can also give an equivalent thermal interpretation of   (\ref{eq:conductivity}) by replacing the word ``voltage" by ``temperature", ``current" by ``heat flux", etc. A mechanical interpretation of   (\ref{eq:conductivity})  is similarly simple: 
it states that $c$ is the point on a string for which the pieces $ [0,c] $ and $[c,1] $ have equal Hooke's constants. By the definition, on the other hand, $c$  is the point which is ``easiest" to slide in the $x$--direction if a spring is grabbed at $c$.

 \vskip 0.3 in 
 \noindent{\bf Remark 2. \,} Operator   (\ref{eq:oned}) arises in many physical settings of which we describe briefly three.
 \begin{enumerate} 
 \item Transversal vibrations  of a string with variable linear density $ \rho (x) $ are governed (to the leading order) by 
 the wave equation 
 $$
 	\rho (x)u_{tt} = u_{xx}.
 $$
 By introducing the mass parameter $s$ via  $  s =\int_{0}^{x}\rho (y)\;dy$  we rewrite the above ODE in the form 
  $$
   u_{tt} = ( a(s) u ^\prime ) ^\prime, \ \ a(s) = \rho (x(s)), \ \ \ ^\prime = \frac{\partial}{\partial s}. 
 $$
\item Longitudinal vibrations of an inhomogeneous string. Consider an elastic string fixed at two ends and undergoing (small) longitudinal vibrations along the $x$--axis, with the ends of the string fixed. Let $ u(x,t) $ denote the displacement from the equilibrium of that particle of the string whose equilibrium position is $x$. Assuming that the material of the string satisfies linear stress-strain relationship, we can rewrite Newton's second law as 
\begin{equation} 
	\rho  \ddot u =(EAu ^\prime ) ^\prime, \ \ \ \dot {} = \frac{\partial}{\partial t},  \ \ \ ^\prime = \frac{\partial}{\partial x} , 
	\label{eq:waveq}
\end{equation}
	where $ E $ is Young's modulus, $\rho= \rho (x) $ is the linear density of the string and $A=A(x) $  is the variable cross-sectional area of the string. By choosing  the new space variable $s$ just as in the preceding example, one reduces   (\ref{eq:waveq})  to the same form as  above:
\begin{equation} 
	  \ddot u =(au ^\prime ) ^\prime. 
	\label{eq:waveq1}
\end{equation}
One can interpret the last equation directly as describing longitudinal vibrations of a string with constant 
linear density $ \rho =1 $ and with variable ``local" Hooke's constant $ a=a(x) $. 
\item The parabolic PDE
\begin{equation} 
	\dot u =(au ^\prime ) ^\prime
 	\label{eq:heateq}
\end{equation}   
	 is the heat equation describing the evolution of temperature $u$ along a rod (with insulated walls) with   heat conductivity $a(x) $ and with heat capacity (specific heat) equal to one unit per unit length\footnote{one can reduce the equation with a variable heat capacity to this one by the same transformation as in the preceding example.}. Indeed,   interpreting $a$ as the heat conductivity amounts to saying that 
$ -a u ^\prime $ is the heat flux along the rod at $x$.  The instantaneous rate of heat gain by a segment $ [x,x+dx] $ is 
then $( -au ^\prime ) ^\prime dx$, to the leading order; since the heat capacity is $1$, this results in   (\ref{eq:heateq}).  
\item The ODE  
$$
	Lu=	(a u^\prime ) ^\prime = 0 
$$
governs an electrostatic potential $u$  along a resistive wire with local conductivity 
$ a(x) $.  Indeed, 
$- a(x) u(x) ^\prime  $ is the current  through point $x$, so that  $ (au ^\prime ) \bigg|_x^{x+dx} =0 $ expresses the conservation of charge in a segment $[x,x+dx] $ (Kirchhoff's first law).   
\end{enumerate}      

 \vskip 0.3 in 
\section{The main result.} 
\label{mr}
In the following theorem we consider the operator $L$ given by   (\ref{eq:oned}). As before, we denote  by $ c=c[L] $   the   point    of least capacity  and by   $ m=m[L]$  the point of maximum of the principal eigenfunction (``the warmest point").  
\begin{thm} \label{thm:main} Consider the elliptic operator $ Lu=	(a u^\prime ) ^\prime $ on $ [0,1]$, where
 $a>0$  is  a monotone increasing   $ C ^{(2)} $--function on the interval.  Then $ c[L]<m[L] $, where 
 $c[L]$ is the point of least capacity associated with $L$, as defined above, and where $x=m[L] $ gives the maximum to the principal eigenfunction of $L$.  
 \end{thm}  
 \begin{figure}[thb]
         \psfrag{c}[c]{\small $c[L]$}
     \psfrag{m}[c]{\small $m[L]$}
     \psfrag{fm}[c]{\small fundamental mode}
     \includegraphics[scale=0.7]{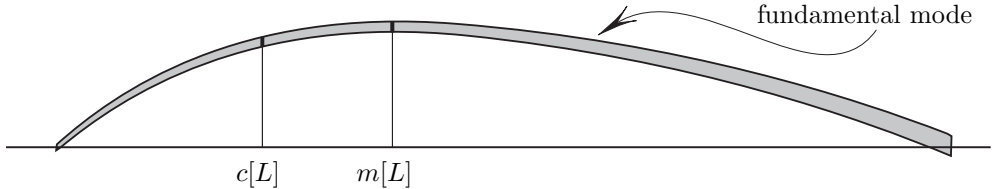}
    \caption{For a tapered vibrating string, the most elastic point is closer to the thin end than  the maximum amplitude point of the fundamental mode.  }
    \label{fig:string}
\end{figure} 

   \begin{cor} {\rm
Consider   vibrations of an inhomogeneous  elastic string, either longitudinal or transversal, described by  $ u_{tt}=
(a(x) u ^\prime ) ^\prime $.   Assume that the ``stiffness coefficient" $ a(x) $    increases from left to right, just as in the last theorem.  Then the the point of maximal amplitude of the fundamental mode lies between the ``softest point " $c$ and the thick end of the string. The same holds for any interval between two consecutive nodes of any overtone.    } 
\end{cor}  
Figure~\ref{fig:string} illustrates this effect for transversal oscillations. 
\vskip 0.3 in 
\noindent{\bf Proof of the theorem.} Let $u (x) $ be the first eigenfunction of the operator (\ref{eq:oned}):
\begin{equation} 
     (a(x) u ^\prime ) ^\prime=- \lambda u, \ \lambda > 0,  \label{eq:eigen}
\end{equation} 
which means that $ u> 0 $ on $ (0,1) $ and satisfies the homogeneous boundary conditions. 
 We rewrite this relation as a system
\begin{equation}  
   \left\{ \begin{array}{l} 
    a(x)u ^\prime = v \\[3pt] 
   \ \ \  \ \  v ^\prime =  - \lambda u  \end{array} \right.   \label{eq:phaseplane}
\end{equation}  
Let us translate the definitions of $m[L]$ and $c[L]$ into geometrical terms. 
Since $c[L]$ is defined by the condition   (\ref{eq:conductivity}), it is natural to  choose 
  $ t=t(x)=  \int_{0}^{x}  a ^{-1} (s) ds $ (the resistance)
as the new independent variable, so that   (\ref{eq:phaseplane})  become 
\begin{equation} 
   \left\{ \begin{array}{l} 
    \dot  U  = V \\[3pt] 
  \dot  V =  - \lambda a(x(t))Y , \end{array} \right.   \label{eq:phaseplane1}
\end{equation} 
where $ \dot{} = \frac{d}{dt} $ and $ U(t)=u(x(t))$, $ V(t)=v(x(t))$.
 Assume without the loss of generality that 
$$
	  \int_{0}^{1}  a ^{-1} (s) ds = 1,   
$$
so that $t$ ranges in $ 0\leq t \leq 1 $. Then $c=c[L]$ is defined by  
 $$
   t(c)= \frac{1}{2}.   
 $$  
 On the other hand, 
 since $x=m$ is the maximizer of $u(x)$,  we have   $ v(m)=u ^\prime (m)=0 $, i.e. 
 $ V(t(m)) =0 $. Introducing the angle
   $\theta = \arg (U+iV)$ we restate the definition of $ x=m $ as $ \theta (t(m))=0 $.    
   Since $ t $ is a monotone increasing function of $x$, proving $ c<m $ is equivalent to proving that
   $ t(c)< t(m) $. Since $ t(c)= \frac{1}{2} $, the proof reduces to   showing that  
   \begin{equation} 
	\frac{1}{2} < t_m.  
	\label{eq:ineq}
\end{equation}   
	
	\begin{figure}[thb]
         \psfrag{1}[c]{\small $ t_m- \tau ^\ast$}
	\psfrag{2}[c]{\small $t_m$}
	\psfrag{3}[c]{\small $\theta (t_m+\tau^\ast)=\theta (1) = \pi /2$}
	\psfrag{4}[c]{\small $\varphi( \tau^ \ast) $}
	\psfrag{5}[c]{\small $\psi( \tau^ \ast) $}
	\psfrag{6}[c]{\small $t= \frac{1}{2} $?}
	\center{  \includegraphics[scale=0.7]{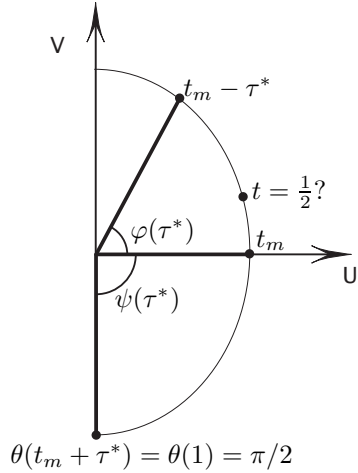}}
	\caption{Proving that $ c[L]<m[L]$, i.e. that $ \frac{1}{2} < t_m $.   }
	\label{fig:phaseplane}
\end{figure} 
Now the angle $\theta$ satisfies 	
 \begin{equation} 
     \dot \theta = - \sin ^2 \theta - \lambda   a(x(t)) \cos ^2 \theta \ 
     \buildrel{def}\over{=}\   f( \theta , t )       
     	\label{eq:argument}
\end{equation} 
 along  with the boundary conditions $     \theta (0)= \frac{\pi}{2} $ and $ \theta (1) =- \frac{\pi}{2}  $. Note that $a(x(t))$ is a monotone increasing function and that $ f< 0 $ for all values of its arguments, so that $ \theta $ is monotone decreasing.  

 The key to   the proof of   (\ref{eq:ineq})   is  the monotonicity of $a(x(t))  $ and the property
 $ f( \theta , t) = f ( - \theta , t) $.  The idea is to compare the angle $\theta$ at equal times $\tau$  before and after crossing the zero value, as in Figure~\ref{fig:phaseplane}. That is, we   introduce
  \begin{equation} 
\varphi ( \tau )= - \theta (  t_m- \tau  ) \ \ \hbox{and} \ \ \ \psi ( \tau ) = \theta (t_m+ \tau ), 
\end{equation} 
so that $ \varphi (0)= \psi(0)=0 $. Moreover, 
these $ \varphi , \ \psi $ satisfy 
\begin{equation}  
   \left\{ \begin{array}{l} 
   \dot \varphi = \dot \theta (t_m- \tau ) = f ( - \varphi , t_m- \tau )=f( \varphi, t_m- \tau ) \\[3pt] 
    \dot \psi =  f( \psi, t_m+\tau ), \end{array} \right.   
\end{equation}  
 where $ f $ has been defined in (\ref{eq:argument}).
 Since  $A(t) =a(x(t) $ is monotone increasing, $f( \theta , t) $ is monotone decreasing in $t$ and   we have 
  $ f( \varphi  , t_m- \tau ) >f ( \varphi  , t_m+ \tau ) $. Since $\varphi $, $ \psi $ share the initial condition,  the comparison theorem applies: 
 \begin{equation} 
    - \frac{\pi }{2} \leq  \psi( \tau )  <\varphi ( \tau ) <0  
    \label{eq:ps}
\end{equation} 
 for  $ \tau > 0 $,  $0<t_m+ \tau \leq 1$. Thus for some $ \tau = \tau ^\ast$  we have 
$$ 
	 - \frac{\pi }{2} = \psi( \tau^\ast )  <\varphi ( \tau^\ast ) <0 
$$
 or, recalling the definition of $ \varphi $ and $ \psi $: 
$$
	- \frac{\pi }{2}\buildrel{A}\over{=}  \theta (t_m+t^\ast) < - \theta (t_m-t^\ast)< 0. 
$$ 
By ($A$)  we have 
\begin{equation} 
	t_m+ \tau ^\ast=1.
	\label{eq:one}
\end{equation}    
 On the other hand, from the above we have $   \theta (t_m- \tau^\ast)< \pi /2 $, which implies   $ t_m-\tau^\ast> 0 $ (since $ \theta (0)= \pi /2 $ and $ \dot \theta < 0 $). Adding this to   (\ref{eq:one})
 we obtain  $ t_m> \frac{1}{2} $, thus completing the proof.

 \section{An open problem.} 
 \label{op}
\begin{figure}[thb]
         \psfrag{g}[c]{\small $w$}
  \center{ \includegraphics[scale=0.7]{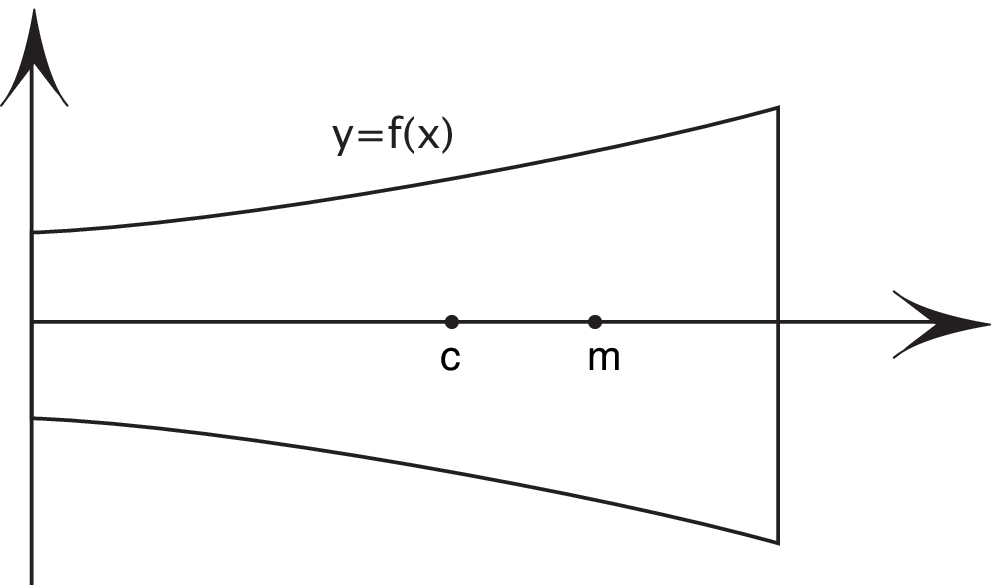}}
    \caption{ $ c<m $?}
    \label{fig:open}
\end{figure}

 Consider a planar domain $D$ (Figure~\ref{fig:open})   bounded by two curves $ y = \pm f(x) $,  where $f$ is a positive monotone increasing function,  and by two lines $ x= 0, \ x=1 $.   
 By symmetry, both points $ c(D) $ and $ m(D) $ lie on the $x$--axis; let us denote by 
   their $x$--coordinates by $c$ and $m$ respectively. 
  
 \vskip 0.3 in 
 \noindent {\it   Open problem 1:}   show that  $ c < m  $.   \vskip 0.3 in 
 \noindent {\it   Open problem 2:}   show that  $ c< m $ for the Laplacian with mixed Dirichlet-Neumann boundary values: 
 Dirichlet on $ x=0$, $ x=1 $ and Neumann on $ y = \pm f(x) $.  
\vskip 0.3 in 
\noindent{\bf Remark.}   If we replace $ y = f(x) $ by $ y = \varepsilon f(x) $ with a small $\varepsilon$, then the  membrane modeled by our   mixed Dirichlet-Neumann problem becomes so narrow as   to resemble an
elastic string described in Theorem \ref{thm:main}. This theorem in fact suggested the second  open problem. For small $\varepsilon$ our  mixed Dirichlet-Neumann problem is approximated by the one--dimensional operator from Theorem \ref{thm:main}. There is a vast literature on partial differential 
operators on thin domains; we  refer to \cite{raugel} and to references therein. 
  
 \newpage
 \bibliography{master}
 \end{document}